\documentclass{article}

\usepackage{makeidx}
\makeindex            

\usepackage{amsmath}  
\usepackage{amssymb}  
\usepackage{dsfont}
\usepackage{amsthm}   
\usepackage{amsfonts}
\usepackage{mathrsfs}

\usepackage{mathrsfs}  
\usepackage{paralist}                       

\newcommand{\depth}[0]{\operatorname{depth}}

\newtheorem{theorem}{Theorem}[section]

\newtheorem*{acknow}{Acknowledgement}

\newtheorem{remark}[theorem]{Remark}

\newcommand*{\hfillplus}{\hfill\linebreak[3]\hspace*{\fill}}

\author{M.~Hellus and R.~H\"ubl}
\title{A result on Macaulay's curve}
\date{\today}   

\begin{document}

\maketitle

\begin{abstract}

We are able to improve what is known about two assumed homogeneous polynomials cutting out Macaulay's curve $C_4\subseteq \mathds P^3_k$ set-theoretically, in characteristic zero. We use local cohomology and an idea from Thoma.

\end{abstract}

\setcounter{section}{-1}

\section{Introduction}

Let $k$ be an algebraically closed field of characteristic zero; let $w,x,y,z$ be indeterminates and
\[ \frak p= \frak p_{C_4}\subseteq k[w,x,y,z]=:R\]
the ideal of Macaulay's curve, that is the curve with parametrization $[s^4:s^3t:st^3:t^4]$ in $\mathds P^3_k$. We assume throughout the paper that $f$ and $g$ are homogeneous polynomials of degrees $d_1$ resp. $d_2$ (in the usual sense) in $\frak p$ that cut out $p$ set-theoretically in the sense that $\sqrt{(f,g)R}= \frak p$. It is a well-known and hard problem to find out whether such polynomials $f$ and $g$ exist or not (such polynomials do exist in positive characteristic: See \cite{H} or \cite[II]{RS}). See \cite{L} for a survey on set-theoretic complete intersections.

\begin{sloppypar}It is natural to endow $R$ with the bigrading $\deg w=(4,0), \deg x=(3,1), \deg y=(1,3), \deg z=(0,4)$. The ideal $p$ in $R$ is bihomogeneous, where ``bihomogeneous'' here and in the following means ``homogeneous with respect to the above bigrading''. We decompose $f$ and $g$ as sums of their bihomogeneous components:\end{sloppypar}
\[ f=\sum_{i=i_\text{min}}^{i_\text{max}}f^{(i,4d_1-i)}, g=\sum_{j=j_\text{min}}^{j_\text{max}}g^{(j,4d_2-j)}\]
($f^{(i,4d_1-i)}$ is the homogeneous component of $f$ of bidegree $(i,4d_1-i)$; similarly for $g$. $i_\text{min}, i_\text{max}, j_\text{min}, j_\text{max}$ are chosen such that $f_\text{min}:=f^{(i_\text{min},4d_1-i_\text{min})}, f_\text{max}:=f^{(i_\text{max},4d_1-i_\text{max})}, g_\text{min}:=g^{(j_\text{min},4d_2-j_\text{min})}, g_\text{max}:=g^{(j_\text{max},4d_2-j_\text{max})}\neq 0$). This situation was studied by Thoma (see e.~g. \cite{T4}, \cite{T3}):

\begin{enumerate}[(i)]

\item Both pairs $f_\text{min},g_\text{min}$ and $f_\text{max},g_\text{max}$ have a proper greatest common divisor in $R$, i.~e. each of the ideals $(f_\text{min},g_\text{min})R$ and $(f_\text{max},g_\text{max})R$ is contained in a prime ideal of height one (\cite[Th.~2.10]{T4}).

\item One of those two greatest common divisors is contained in $\frak p$ (\cite[Cor.~2.9]{T4}).

\end{enumerate}

We are able to improve these results:

\begin{enumerate}[(i)]

\item[(i')] For both ideals $(f_\text{min},g_\text{min})R$ and $(f_\text{max},g_\text{max})R$ it is true that {\it all} minimal prime divisors have height one, i.~e. both ideals are principal up to radical (theorem \ref{theo_1}.(a)).

\item[(ii')] The greatest common divisors of {\it both} pairs $f_\text{min},g_\text{min}$ and $f_\text{max},g_\text{max}$ are contained in $\frak p$ (theorem \ref{theo_1}.(b)).

\end{enumerate}

For our proofs of (i') and (ii') we use local cohomology and a modification of the following idea from Thoma (\cite{T4}):

Let $\lambda$ and $\mu$ be additional indeterminates.
\[ {F(w,x,y,z,\lambda,\mu):=f(\lambda^4w,\lambda^3\mu x,\lambda\mu^3 y,\mu^4z)=\sum_if^{(i,4d_1-i)}\lambda^i\mu^{4d_1-i}\atop G(w,x,y,z,\lambda,\mu):=g(\lambda^4w,\lambda^3\mu x,\lambda\mu^3 y,\mu^4z)=\sum_ig^{(i,4d_2-i)}\lambda^i\mu^{4d_2-i}}\]
$\in R[\lambda,\mu]=:S$. There can be no point $[w_0:x_0:y_0:z_0]\in \mathds P^3_k\setminus C_4$ with $\lambda_0,\mu_0 \in k^*$ and such that
\[ \underbrace{F(w_0,x_0,y_0,z_0,\lambda_0,\mu_0)}_{=f(\lambda_0^4w_0,\lambda_0^3\mu_0x_0,\lambda_0\mu_0^3y_0,\mu_0^4z_0)}=\underbrace{G(w_0,x_0,y_0,z_0,\lambda_0,\mu_0)}_{=g(\lambda_0^4w_0,\lambda_0^3\mu_0x_0,\lambda_0\mu_0^3y_0,\mu_0^4z_0)}=0,\]
i.~e.
\[ [\lambda_0^4w_0:\lambda_0^3\mu_0x_0:\lambda_0\mu_0^3y_0:\mu_0^4z_0]\]
belongs to $V(f,g)$, because such a point could not have the form $[s^4:s^3t:st^3:t^4]$ (otherwise -- because of $\lambda_0,\mu_0\neq 0$ -- also $[w_0:x_0:y_0:z_0]$ would have this form and would therefore belong to $C_4$).

\begin{acknow}

We thank Peter Schenzel for a comment which lead to a substantial simplification in the proof of theorem \ref{theo_1}.

\end{acknow}

\section{Results}

We modify Thoma's observation which was described in the introduction: Recall that we assume that $f$ and $g$ are homogeneous polynimals of degrees $d_1$ resp. $d_2$ in $\frak p= \frak p_{C_4}$ that cut out $p$ set-theoretically in the sense that $\sqrt{(f,g)R}= \frak p$. Let $\lambda$ be an additional indeterminate.
\[ {F_1(w,x,y,z,\lambda):=f(\lambda^4w,\lambda^3x,\lambda y,z)=\sum_if^{(i,4d_1-i)}\lambda^i\atop G_1(w,x,y,z,\lambda):=g(\lambda^4w,\lambda^3x,\lambda y,z)=\sum_ig^{(i,4d_2-i)}\lambda^i}\]
$\in R[\lambda]=: T$.  There can be no point $[w_0:x_0:y_0:z_0]\in \mathds P^3_k\setminus C_4$, represented by  $(w_0,x_0,y_0,z_0)\neq (0,0,0,0)$, and no $\lambda_0\in k^*$  such that
\[ \underbrace{F_1(w_0,x_0,y_0,z_0,\lambda_0)}_{=f(\lambda_0^4w_0,\lambda_0^3 x_0,\lambda_0 y_0,z_0)}=\underbrace{G_1(w_0,x_0,y_0,z_0,\lambda_0)}_{=g(\lambda_0^4w_0,\lambda_0^3 x_0,\lambda_0 y_0,z_0)}=0,\]
i.~e.
\[ [\lambda_0^4w_0:\lambda_0^3 x_0:\lambda_0 y_0:z_0]\]
belongs to $V(f,g)$, because such a point could not have the form $[s^4:s^3t:st^3:t^4]$ (otherwise -- because of $\lambda_0\neq 0$ -- also $[w_0:x_0:y_0:z_0]$ would have this form and would therefore belong to $C_4$).

\begin{remark}

\label{dual_evr}All arguments in the sequel can be translated in an obvious way to
\[ {F_2(w,x,y,z,\mu):=f(w,\mu x,\mu^3 y,\mu^4 z)=\sum_if^{(i,4d_1-i)}\mu^{4d_1-i}\atop G_2(w,x,y,z,\mu):=g(w,\mu x,\mu^3 y,\mu^4z)=\sum_ig^{(i,4d_2-i)}\mu^{4d_2-i}}\]
$\in R[\mu]$ (this time, of course, $\mu$ is the additional indeterminate); the obtained results are analogous with 'max' instead of 'min'.

\end{remark}

\hfillplus $\square$

Note that the observation from the beginning of this section works equally for $\tilde F_1$ and $\tilde G_1$, where these two polynomials are obtained from $F_1$ and $G_1$ by cancelling out $\lambda$ as much as possible. We claim that
\begin{equation}\label{mpd1} \sqrt {(\tilde F_1,\tilde G_1)T}=\sqrt{(\lambda, f_\text{min},g_\text{min})T}\cap \frak p_{C_4}T\end{equation}
(with $T = R[\lambda]$ as above).

``$\subseteq $'' is trivial; ``$\supseteq $'': Let $P=(w_0,x_0,y_0,z_0,\lambda_0)$ be an arbitrary point on $V(\tilde F_1,\tilde G_1)$, we have to show that $P\in V(\lambda, f_\text{min},g_\text{min})\cup V(\frak p_{C_4}T)$: In case $\lambda_0=0$, evaluation at $P$ makes all bihomogeneous components of $\tilde F_1$ and of $\tilde G_1$ apart from ``the minimal ones'' vanish, therefore we have $P\in V(\lambda, f_\text{min},g_\text{min})$; and the other case $\lambda_0\neq 0$ follows from the observation from the beginning of this section.

We study the minimal prime divisors of $(\tilde F_1,\tilde G_1)T$, our main source of information is formula (\ref{mpd1}).

Let
\[ \frak p_1,\ldots ,\frak p_r, \frak q_1,\ldots, \frak q_s\]
\begin{sloppypar}be exactly the minimal prime divisors of $(f_\text{min},g_\text{min})R$, where all $\frak p_i$ have height one (these $\frak p_i$ are therefore principal, they encode information on the $\gcd (f_\text{min},g_\text{min})$) and all $\frak q_i$ have height two. Clearly,\end{sloppypar}
\[ (\lambda, \frak p_1),\ldots ,(\lambda, \frak p_r),(\lambda, \frak q_1),\ldots ,(\lambda, \frak q_s)\]
are exactly the minimal prime divisors of $(\lambda, f_\text{min},g_\text{min})T$.

The prime ideals $(\lambda, \frak p_i)$ have height two and the prime ideals $(\lambda, \frak q_i)$ have height three. We get
\begin{equation}\label{mpd2}\sqrt{(\tilde F_1,\tilde G_1)T}=[(\lambda, \frak p_1)\cap \ldots \cap (\lambda, \frak p_r)\cap (\lambda,\frak q_1)\cap \ldots \cap(\lambda,\frak q_s)]\cap \frak p_{C_4}T\end{equation}
All prime ideals occurring in the bracket $[\ldots ]$ contain $\lambda $, $\frak p_{C_4}T$ does not contain $\lambda$.

In particular, between these $r+s+1$ prime ideals only one type of inclusion is possible: $\frak p_{C_4}T$ can possibly be contained in some $(\lambda,\frak q_i)$, equivalently: $\frak p_{C_4}$ is contained in (and therefore equal to) some $\frak q_i$.

\begin{itemize}

\item First case: $\frak p_{C_4}$ is contained in no $\frak q_i$: This means that between our $r+s+1$ prime ideals there are no inclusions at all. Since there are no inclusions, formula (\ref{mpd2}) is the (unique) minimal decomposition of $\sqrt{(\tilde F_1,\tilde G_1)T}$, the $r+s+1$ prime ideals are exactly the minimal prime divisors of $(\tilde F_1,\tilde G_1)T$. But the latter ideal has no minimal prime divisors of height three, i.~e. $s=0$. All minimal prime divisors of $(f_\text{min},g_\text{min})R$ have height one.

\item Second case: $\frak p_{C_4}$ is contained in (and therefore equal to) some $\frak q_i$: There is exactly one inclusion among the prime ideals in (\ref{mpd2}), namely $\frak p_{C_4}\subseteq (\lambda,\frak q_i)$; since this is the only inclusion, omitting $(\lambda,\frak q_i)$ from (\ref{mpd2}) leads to the minimal decomposition of $\sqrt{(\tilde F_1,\tilde G_1)T}$. However, again, no minimal prime divisor of $(\tilde F_1,\tilde G_1)T$ has height three. Therefore, we must have $s=1$, i.~e. $(f_\text{min},g_\text{min})R$ has exactly one minimal prime divisor of height two, namely $\frak p_{C_4}$ (we will see below that this second case is actually impossible).

\end{itemize}

\begin{sloppypar}We summarize: The only minimal prime divisor of height two which $(f_\text{min},g_\text{min})R$ can possibly have, is $\frak p_{C_4}$. Furthermore, at least one minimal prime divisor of height one must exist, since otherwise the radical of $(f_\text{min},g_\text{min})R$ would equal its (only) minimal prime divisor $\frak p_{C_4}$, contradicting \cite[p.~816]{T3}. In particular, we may write (\ref{mpd1}) in the form\end{sloppypar}

\begin{equation}\label{mpd3} \sqrt {(\tilde F_1,\tilde G_1)T}=\sqrt{(\lambda, t_\text{min})T}\cap \frak p_{C_4}T.\end{equation}

with $t_\text{min}$ the greatest common divisor of $f_\text{min}$ and $g_\text{min}$ ($t_\text{min}$ is not a unit since $(f_\text{min},g_\text{min})R$ has a minimal prime divisor of height one).

\begin{theorem}

\label{theo_1}\begin{enumerate}[(a)]

\item \begin{sloppypar}All minimal prime divisors of $(f_\text{min},g_\text{min})R$ and of $(f_\text{max},g_\text{max})R$ have height one. In particular, both pairs $f_\text{min},g_\text{min}$ and $f_\text{max},g_\text{max}$ have a proper (non-unit) common divisor $t_\text{min}$ resp. $t_\text{max}$ in $R$.\end{sloppypar}

\item Both $t_\text{min}$ and $t_\text{max}$ are contained in $\frak p_{C_4}$.

\item \begin{equation}\label{mpd4} \sqrt {(\tilde F,\tilde G)S}=(\lambda, \mu)\cap \sqrt{(\lambda, t_\text{min})S}\cap \sqrt{(\mu, t_\text{max})S}\cap \frak p_{C_4}S.\end{equation}

\end{enumerate}

\end{theorem}

{\it Proof:} (a) and (b): The ring $S_1=k[\lambda]_{(\lambda)}[w,x,y,z]=R_0[w,x,y,z]$ with $R_0$ the subring of elements of degree zero and $\deg w,x,y,z=1$ is graded and *local (in particular: noetherian), using terminology from \cite{BS}. It is also a localization of $T$. The ring $\overline{S_1}:=S_1/(\tilde F_1, \tilde G_1)S_1$ is also *local and we can formulate (\ref{mpd3}) for its ideals $a:=\sqrt{(\lambda, t_\text{min})}\overline{S_1}$ and $b:=\frak p_{C_4}\overline{S_1}$. (\ref{mpd3}) says that $ab$ is nilpotent, however $a$ and $b$ are non-nilpotent (this is clear from the discussion preceeding (\ref{mpd3})). The following trick is known, to the best of our knowledge it goes back to Irving Kaplansky (see also \cite[Prop.~2.1]{H}):
The (exact) Mayer-Vietoris sequence for local cohomology of $\overline{S_1}$ with respect to $a$ and $b$ starts as follows:
\[ \tag{$*$}0\to \Gamma_{a+b}(\overline{S_1})\to \Gamma_a(\overline{S_1})\oplus \Gamma_b(\overline{S_1})\to \Gamma_{a\cap b}\overline{S_1}=\overline{S_1}\to H^1_{a+b}(\overline{S_1}).\]
Therefore, $\depth (a+b,\overline{S_1})$ must be at most one, because otherwise ($*$) would provide an isomorphism
\[ \Gamma_a(\overline{S_1})\oplus \Gamma_b(\overline{S_1})\cong \overline{S_1},\]
which is impossible since $\overline{S_1}$ is *local. In the ring $S_1$ this means that the depth and hence also the height of $\sqrt{(\lambda, t_\text{min})}S_1+\frak p_{C_4}S_1$ is at most $1+2=3$ (note that $\tilde F_1, \tilde G_1$ is a regular sequence in $S_1$, e.~g. by \ref{mpd3}). But this is only possibly if $t_\text{min}$ is in $p_{C_4}$.

Now, both $(f_\text{min},g_\text{min})R$ and $(f_\text{max},g_\text{max})R$ have a minimal prime divisor of height one and which is contained in $\frak p_{C_4}$; therefore $\frak p_{C_4}$, which is -- as we have seen above -- the only possible minimal prime divisor of height two, cannot be minimal. Consequently, all minimal prime divisors of $(f_\text{min},g_\text{min})R$ (analougously, of $(f_\text{max},g_\text{max})R$) have height one.

{\it (c)} We work with the polynomials $F$ and $G$ from the introduction. The observation from the end of the introduction works equally for $\tilde F$ and $\tilde G$, where these two polynomials are obtained from $F$ and $G$ by cancelling out $\lambda$ and $\mu$ as much as possible. We claim that
\begin{equation}\label{mpd5} \sqrt {(\tilde F,\tilde G)S}=(\lambda, \mu)\cap \sqrt{(\lambda, f_\text{min},g_\text{min})S}\cap \sqrt{(\mu, f_\text{max},g_\text{max})S}\cap \frak p_{C_4}S.\end{equation}
\begin{sloppypar}``$\subseteq $'' follows from the fact that both $f$ and $g$ consist of at least two bihomogeneous components (\cite[p.~816]{T3}, \cite[Lemma~3.1]{T3}); ``$\supseteq $'': Let $P=(w_0,x_0,y_0,z_0,\lambda_0,\mu_0)$ be an arbitrary point on $V(\tilde F,\tilde G)$, we have to show that $P\in V(\lambda,\mu)\cup V(\lambda, f_\text{min},g_\text{min})\cup V(\mu, f_\text{max},g_\text{max})\cup V(\frak p_{C_4}S)$: The case $\lambda_0,\mu_0=0$ is trivial; case $\lambda_0=0, \mu_0\neq 0$: Evaluation at $P$ makes all bihomogeneous components of $\tilde F$ and of $\tilde G$ apart from ``the minimal ones'' vanish, therefore we have $P\in V(\lambda, f_\text{min},g_\text{min})$; the case $\mu_0=0,\lambda_0\neq 0$ is analogous with ``maximal'' instead of ``minimal''; finally the case $\lambda_0\neq 0,\mu_0\neq 0$ follows from Thoma's idea described in the introduction.\end{sloppypar}\hfillplus $\square$

\begin{remark}

It is clear that $\lambda\cdot t_\text{max}$ (and, similarly, $\mu\cdot t_\text{min}$) belongs to all four ideals in the right hand side of formula (\ref{mpd3}) and, therefore, a power of it can be written as a linear combination of $\tilde F$ and $\tilde G$. Similarly, $t_\text{min}$ is in the radical of $(f_\text{min},g_\text{min})$ and  $t_\text{max}$ is in the radical of $(f_\text{max},g_\text{max})$.

\end{remark}

\begin{remark}

The non-trivial result from Thoma that the number of (non-zero)
bihomogenous components of $f$ or of $g$ is at least three
(\cite[Th.~3.1.(a)]{T4}) follows immdediately from theorem \ref{theo_1} b)
together with the well-known fact that the number of (non-zero) bihomogenous
components of $f$ or of $g$ is at least two (\cite[p.~816]{T3}). (Is also
well-known and easy to see that neither $f$ nor $g$ can be bihomogenous
(\cite[Lemma~3.1]{T3})).

\end{remark}

\begin{remark}

Finally note that all our results and proofs immediately generalize to arbitrary symmetric, non-arithmetically Cohen Macualay monomial curves $[s^d:s^at^b:s^bt^a:t^d]$ (see \cite[p.~816]{T3}).

\end{remark}


\begin{thebibliography}{1}

\bibitem[BS]{BS} Brodmann, Markus P., and Rodney Y. Sharp. Local cohomology: an algebraic introduction with geometric applications. Vol. {\bf 136}. Cambridge University Press, 2012.

\bibitem[H]{H} R. Hartshorne, Complete intersections in characteristic $p>0$, Amer. J. Math. {\bf 101} (1979), no.~2, 380--383.

\bibitem[L]{L} G. Lyubeznik, A survey of problems and results on the number of defining equations, in {\it Commutative algebra (Berkeley, CA, 1987)}, 375--390, Math. Sci. Res. Inst. Publ., 15 Springer, New York.

\bibitem[RS]{RS} H. Roloff\ and\ J. St\"uckrad, Bemerkungen \"uber Zusammenhangseigenschaften und mengentheoretische Darstellung projektiver algebraischer Mannigfaltigkeiten, Wiss. Beitr. Martin-Luther-Univ. Halle-Wittenberg M {\bf 12} (1979), 125--131.

\bibitem[T1]{T1} A. Thoma, Monomial space curves in ${\bf P}\sp 3\sb k$ as binomial set theoretic complete intersections, Proc. Amer. Math. Soc. {\bf 107} (1989), no.~1, 55--61.

\bibitem[T2]{T2} A. Thoma, On set-theoretic complete intersections in ${\bf P}\sp 3\sb k$, Manuscripta Math. {\bf 70} (1991), no.~3, 261--266.

\bibitem[T3]{T3} A. Thoma, On the arithmetically Cohen-Macaulay property for monomial curves, Comm. Algebra {\bf 22} (1994), no.~3, 805--821.

\bibitem[T4]{T4}A. Thoma, On the equations defining monomial curves, Comm. Algebra {\bf 22} (1994), no.~7, 2639--2649.

\end{thebibliography}
\end{document}